\documentclass[10pt]{amsart}
\usepackage[leqno]{amsmath}
\usepackage{amssymb,latexsym,soul,cite,amsthm,color,enumitem,graphicx,mathtools,microtype,accents}
\usepackage[colorlinks=true,urlcolor=airforceblue,citecolor=airforceblue,linkcolor=airforceblue,linktocpage,pdfpagelabels,bookmarksnumbered,bookmarksopen]{hyperref}
\definecolor{airforceblue}{rgb}{0.36, 0.54, 0.66}
\usepackage[english]{babel}
\usepackage[left=2.5cm,right=2.5cm,top=2.5cm,bottom=2.5cm]{geometry}

\numberwithin{equation}{section}
\newtheorem{theorem}{Theorem}[section]
\theoremstyle{plain}
\newtheorem{lemma}[theorem]{Lemma}
\theoremstyle{plain}
\newtheorem{proposition}[theorem]{Proposition}
\theoremstyle{plain}

\theoremstyle{definition}
\newtheorem{remark}[theorem]{Remark}

\newtheorem{problem}[theorem]{Problem}

\newcommand{\N}{{\mathbb N}}

\newcommand{\R}{{\mathbb R}}
\newcommand{\eps}{\varepsilon}
\newcommand{\beq}{\begin{equation}}
\newcommand{\eeq}{\end{equation}}
\renewcommand{\le}{\leqslant}
\renewcommand{\ge}{\geqslant}

\def\Xint#1{\mathchoice
{\XXint\displaystyle\textstyle{#1}}%
{\XXint\textstyle\scriptstyle{#1}}%
{\XXint\scriptstyle\scriptscriptstyle{#1}}%
{\XXint\scriptscriptstyle\scriptscriptstyle{#1}}%
\!\int}
\def\XXint#1#2#3{{\setbox0=\hbox{$#1{#2#3}{\int}$ }
\vcenter{\hbox{$#2#3$ }}\kern-.6\wd0}}
\def\dashint{\Xint-}

\makeatletter
\newcommand{\leqnomode}{\tagsleft@true}
\newcommand{\reqnomode}{\tagsleft@false}
\makeatother

\newenvironment{enumroman}{\begin{enumerate}

}{\end{enumerate}}

\title[Boundary regularity for fractional $p$-Laplacian]{A survey on boundary regularity for the fractional $p$-Laplacian and its applications / Un'indagine sulla regolarit\`a alla frontiera per il $p$-laplaciano frazionario e le sue applicazioni}

\author[A.\ Iannizzotto]{Antonio Iannizzotto}

\address[]{Dipartimento di Matematica e Informatica
\newline\indent
Universit\`a degli Studi di Cagliari
\newline\indent
Via Ospedale 72, 09124 Cagliari, Italy}
\email{antonio.iannizzotto@unica.it}


\begin{document}

\begin{abstract}
We survey some recent regularity results for fractional $p$-Laplacian elliptic equations, especially focusing on pure and weighted boundary H\"older continuity of the solutions of related Dirichlet problems. Then, we present some applications of such results to general nonlinear elliptic equations of fractional order, treated via either variational or topological methods. / Esaminiamo alcuni recenti risultati di regolarit\`a per equazioni ellittiche con $p$-laplaciano frazionario, concentrandoci specialmente sulla continuit\`a h\"olderiana alla frontiera delle soluzioni dei relativi problemi di Dirichlet. Quindi presentiamo alcune applicazioni di tali risultati a equazioni ellittiche non lineari di ordine frazionario pi\`u generali, trattate con metodi sia variazionali che topologici.
\end{abstract}

\maketitle

\begin{center}
Version of \today\
\end{center}

\section{Introduction and preliminaries}\label{sec1}

\noindent
The purpose of the present note is to survey some recent results concerning regularity of the solutions of elliptic partial differential equations belonging to the following class:
\beq\label{lps}
\mathcal{L}_{p,K}^s u = f(x) \quad \text{in $\Omega$.}
\eeq
Here $\Omega\subset\R^N$ ($N\ge 2$) is an open set, $f\in L^\infty(\Omega)$, and the leading operator is a nonlinear one with fractional order $2s$, admitting the following representation for $u:\R^N\to\R$ smooth enough and $x\in\R^N$:
\[\mathcal{L}_{p,K}^s u(x) = \lim_{\eps\to 0^+}\int_{B_\eps(x)^c}|u(x)-u(y)|^{p-2}(u(x)-u(y))K(x,y)\,dy,\]
where $s\in(0,1)$, $p>1$, and the kernel $K:\R^N\times\R^N\to\R$ is a measurable function s.t.\ for a.e.\ $x,y\in\R^N$
\begin{itemize}[leftmargin=1cm]
\item[$(K_1)$] $K(x,y)=K(y,x)$;
\item[$(K_2)$] $\Lambda_1\le K(x,y)|x-y|^{N+ps}\le\Lambda_2$ ($0<\Lambda_1\le\Lambda_2$).
\end{itemize}
If $\Lambda_1=\Lambda_2=C_{N,p,s}>0$ (a normalization constant which varies from one reference to another), $\mathcal{L}_{p,K}^s$ becomes the $s$-fractional $p$-Laplacian, namely,
\[(-\Delta)_p^s u(x) = C_{N,p,s}\lim_{\eps\to 0^+}\int_{B_\eps(x)^c}\frac{|u(x)-u(y)|^{p-2}(u(x)-u(y))}{|x-y|^{N+ps}}\,dy.\]
If, in addition, $p=2$, then we retrieve the classical fractional Laplacian $(-\Delta)^s$. In the linear case, $p=2$, elliptic equations of the type \eqref{lps} and their evolutive counterparts have countless applications to quantum mechanics, flame propagation, dislocation of crystals, and above all models based on stable L\'evy-type stochastic processes, see \cite{BV} for an elementary introduction to this vast subject.
\vskip2pt
\noindent
In the nonlinear framework, the main motivation arises from game theory, which may lead to either fully nonlinear equations of fractional order, or problems driven by the fractional $\infty$-Laplacian (the latter is also related to the problem of H\"older continuous extensions), see \cite{BCF,C}. Equation \eqref{lps}, however, does not fall in this class, being rather a general divergence-form nonlinear nonlocal equation, of degenerate ($p>2$) or singular ($p<2$) type. The model operator $(-\Delta)_p^s$ can be seen as both an approximation of the $p$-Laplacian for $p$ fixed and $s\to 1$ (see \cite{IN}, and \cite{BBM} for a functional-analytic approach), and an approximation of the fractional infinity Laplacian for $s$ fixed and $p\to\infty$ (see \cite{CLM,LL}).
\vskip2pt
\noindent
Anyway, also due to its independent mathematical interest, the fractional $p$-Laplacian has become a major subject of research in the last decade. Several results on existence, multiplicity, and qualitative properties have appeared, beginning with the very simple Morse-theoretic approach of \cite{ILPS}, while the regularity theory for such operator has undergone a slower development and is growing fast in the last years. The reason for such a gap is easily understood.
\vskip2pt
\noindent
Indeed, the functional analytic properties of $(-\Delta)_p^s$ like continuity, monotonicity, spectral properties, and weak minimum/comparison principles do not differ too much from those of the local $p$-Laplacian, and the same holds for the more general operator $\mathcal{L}_{p,K}^s$, see \cite{ILPS,LL}. On the other hand, regularity theory for fractional order operators requires new ideas with respect to the classical framework, mainly due to two distinctive features: first, the operator does not involve derivatives, but rather a fractional order difference quotient of singular nature; second, the operator has a {\em nonlocal nature}, meaning that perturbing the solution outside the domain $\Omega$ does affect the equation in $\Omega$. The latter, in particular, will be a crucial element in the forthcoming discussion. Some years ago, the state of the art on nonlinear nonlocal equations was portrayed in \cite{MS,P}. We now aim at providing new details, including some very recent results, in a smooth exposition. Here we do not cover the field of evolutive fractional $p$-Laplacian equation, referring the reader to \cite{L1} and the references therein.
\vskip2pt
\noindent
Therefore, in the present note we will focus on some regularity results for solutions of equations of the class \eqref{lps}, especially on the issue of global (or boundary) regularity for the solutions of the Dirichlet problem
\beq\label{dir}
\begin{cases}
\mathcal{L}_{p,K}^s u = f(x) & \text{in $\Omega$} \\
u = 0 & \text{in $\Omega^c$.}
\end{cases}
\eeq
For the model case $K(x,y)=|x-y|^{-N-ps}$, optimal boundary regularity of $u$, which amounts to $s$-H\"older continuity in $\overline\Omega$, was recently proved in \cite{IM}, along with H\"older continuity of the quotient $u/{\rm dist}^s(\cdot,\Omega^c)$. Such study is motivated by some applications to more general nonlinear equations, which we will discuss at the end of the paper.
\vskip2pt
\noindent
First, we need more precise definitions of the involved operators and equations. The pointwise representation given above for $\mathcal{L}_{p,K}^s$ is a bit na\"{\i}ve, since in general (especially for $p<2$) the integral may fail to converge. So, we need to define our operator more precisely in the natural framework of {\em fractional Sobolev spaces}, see \cite{DNPV,L}.
\vskip2pt
\noindent
For any open set $\Omega\subseteq\R^N$, $p>1$, $s\in(0,1)$ we define the Gagliardo seminorm of a measurable function $u:\Omega\to\R$ as
\[[u]_{s,p,\Omega} = \Big[\iint_{\Omega\times\Omega}\frac{|u(x)-u(y)|^p}{|x-y|^{N+ps}}\,dx\,dy\Big]^\frac{1}{p}.\]
We define the Sobolev space
\[W^{s,p}(\Omega) = \big\{u\in L^p(\Omega):\,[u]_{s,p,\Omega}<\infty\big\},\]
a uniformly convex, separable Banach space endowed with the norm
\[\|u\|_{W^{s,p}(\Omega)} = \big(\|u\|_{L^p(\Omega)}^p+[u]_{s,p,\Omega}^p\big)^\frac{1}{p}.\]
If $\Omega$ is bounded, we define the localized space
\[\widetilde{W}^{s,p}(\Omega) = \Big\{u\in L^p_{\rm loc}(\R^N):\,u\in W^{s,p}(\Omega') \ \text{for some $\Omega'\Supset\Omega$,} \ \int_{\R^N}\frac{|u(x)|^{p-1}}{(1+|x|)^{N+ps}}\,dx<\infty\Big\},\]
while for unbounded $\Omega$ we set
\[\widetilde{W}^{s,p}_{\rm loc}(\Omega) = \big\{u\in L^p_{\rm loc}(\R^N):\,u\in W^{s,p}(\Omega') \ \text{for all $\Omega'\Subset\Omega$}\big\}.\]
We say that $u\in\widetilde{W}^{s,p}_{\rm loc}(\Omega)$ is a (local weak) solution of \eqref{lps} if for all $\varphi\in C^\infty_c(\Omega)$
\beq\label{wls}
\iint_{\R^N\times\R^N}|u(x)-u(y)|^{p-2}(u(x)-u(y))(\varphi(x)-\varphi(y))K(x,y)\,dx\,dy = \int_\Omega f(x)\varphi(x)\,dx.
\eeq
Alternatively, weak solutions can be defined by means of the tail spaces introduced in \cite{KKP}, namely
\[L^{p-1}_{ps}(\R^N) = \Big\{u\in L^{p-1}_{\rm loc}(\R^N):\,\int_{\R^N}\frac{|u(x)|^{p-1}}{(1+|x|)^{N+ps}}\,dx<\infty\Big\}.\]
We note that $W^{s,p}(\R^N)\subset L^{p-1}_{ps}(\R^N)$, while for a general domain $\Omega$ such inclusion is not granted. Tail spaces, introduced to deal with obstacle problems, provide a more general framework but at the price of a slightly more complex definition of solutions, so for the purposes of the present survey we will keep the definition based on $\widetilde{W}^{s,p}(\Omega)$. If $\Omega$ is bounded with a $C^{1,1}$-boundary, then the subspace
\[W^{s,p}_0(\Omega) = \big\{u\in W^{s,p}(\R^N):\,u=0 \ \text{in $\Omega^c$}\big\}\]
can be endowed with the norm $\|u\|_{W^{s,p}_0(\Omega)}=[u]_{s,p,\R^N}$ and is compactly embedded into $L^q(\Omega)$, for all $q\in[1,p^*_s)$, where
\[p^*_s = \begin{cases}
\displaystyle\frac{Np}{N-ps} & \text{if $ps<N$} \\
\infty & \text{if $ps\ge N$.}
\end{cases}\]
In such case, $C^\infty_c(\Omega)$ is a dense subspace of $W^{s,p}_0(\Omega)$. We say that $u\in W^{s,p}_0(\Omega)$ is a (weak) solution of \eqref{dir} if $u$ satisfies \eqref{wls} for all $\varphi\in C^\infty_c(\Omega)$. We are going to investigate the regularity of such solutions.

\section{Interior regularity}\label{sec2}

\noindent
In the linear case $p=2$, interior regularity of the solutions to \eqref{lps}-type equations is a well established subject. For the model operator $(-\Delta)^s$, we have the following:

\begin{theorem}\label{irl}
{\rm (Interior regularity, linear case)} Let $\Omega\subset\R^N$ be an open set, $s\in(0,1)$ $f\in L^\infty(\Omega)$, and $u$ be a solution of $(-\Delta)^su=f$ in $\Omega$. Then:
\begin{enumroman}
\item\label{irl1} {\rm (H\"older continuity)} if $s\neq 1/2$ then $u\in C^{2s}_{\rm loc}(\Omega)$, and if $s=1/2$ then $u\in C^\alpha_{\rm loc}(\Omega)$ for all $\alpha\in(0,2s)$;
\item\label{irl2} {\rm (Schauder estimates)} if $f\in C^\alpha(\Omega)$ and $2s+\alpha\notin\N$ then $u\in C^{2s+\alpha}_{\rm loc}(\Omega)$.
\end{enumroman}
\end{theorem}

\noindent
In \ref{irl2} it is understood that, if $2s+\alpha>1$, then $u\in C^1_{\rm loc}(\Omega)$ and $\nabla u\in C^{2s+\alpha-1}_{\rm loc}(\Omega)$. Theorem \ref{irl} follows from \cite[Theorem 2.4.1, Proposition 2.4.4]{FRRO}, see also \cite{BCI,DROSV,FRO,ROS2} for more general linear nonlocal operators invariant by translation,  \cite{C1,DK} for measure-type kernels, and \cite{M} for a discussion on optimal H\"older exponent in \ref{irl1}. Most regularity results come with a uniform estimate on the H\"older norm of solutions. A valuable tool for linear nonlocal equations is the extension result of \cite{CS}, which also allows to prove Harnack inequalities.
\vskip2pt
\noindent
In the nonlinear framework, such tools are not fully developed yet. Recently, in \cite{DTGCV} an extension operator relating the fractional $p$-Laplacian to a local elliptic equation in dimension $N+1$ was introduced, under the condition $p(2-s)<2$ and restricted to $C^2$-functions. Such result, though interesting as a nonlinear counterpart of \cite{CS}, does not provide an equally general and useful tool for regularity theory.
\vskip2pt
\noindent
The first regularity results for nonlinear, nonlocal operators of the form $\mathcal{L}_{p,K}^s$ are found in \cite{DCKP,DCKP1}, where the authors use a nonlocal adaptation of the De Giorgi-Nash-Moser method to prove that, if $f=0$, weak solutions to the equation \eqref{lps} are locally bounded, locally H\"older continuous with an undetermined small exponent, and satisfy a Harnack inequality. A similar result is proved with a different method in the recent paper \cite{CDI}, based on the clustering and expansion of positivity (see also \cite{DIV}). Both approaches involve careful estimates of the {\em nonlocal tail}, a special quantity depending on the behavior of the solution $u$ outside a ball $B_R(x)$, defined by
\[{\rm Tail}(u,x,R) = \Big[R^{ps}\int_{B_R^c(x)}\frac{|u(x)|^{p-1}}{|x-y|^{N+ps}}\,dy\Big]^\frac{1}{p-1}.\]
The control of tail terms is one of the most delicate issues in nonlocal regularity theory. During the last decade, such theory has widely developed. For instance, we mention the H\"older regularity for solutions of non-homogenous equations with measure-type data from \cite{KMS} and a very general result for fractional De Giorgi classes from \cite{C2}. For the case of the fractional $p$-Laplacian, a basic reference result is the following:

\begin{theorem}\label{irn}
{\rm (Interior regularity, nonlinear case)} Let $\Omega$ be a bounded open set, $p>1$, $s\in(0,1)$, $f\in L^\infty(\Omega)$, and $u\in\widetilde{W}^{s,p}(\Omega)$ be a solution of $(-\Delta)_p^s u=f$ in $\Omega$. Then, for all $\alpha$ satisfying
\[0 < \alpha < \min\Big\{1,\,\frac{ps}{p-1}\Big\}\]
we have $u\in C^\alpha_{\rm loc}(\Omega)$. In addition, there exists $C=C(N,p,s,\alpha)>0$ s.t.\ for all $B_{4R}(x)\subset\Omega$ we have
\[[u]_{C^\alpha(B_{R/8}(x))} \le \frac{C}{R^\alpha}\Big[\|u\|_{L^\infty(B_R(x))}+R^\frac{ps}{p-1}\|f\|_{L^\infty(\Omega)}^\frac{1}{p-1}+{\rm Tail}(u,x,R)\Big].\]
\end{theorem}

\noindent
Theorem \ref{irn} was proved in \cite{BLS} for the degenerate case $p>2$ and in \cite{GL} for the singular case $p<2$, respectively. We remark that the results stated above are in fact more general, as they hold for $f\in L^q(\Omega)$, with
\[q > \max\Big\{1,\,\frac{N}{ps}\Big\}\]
and a suitably adapted H\"older exponent $\alpha$ depending on $q$. The H\"older exponent $\alpha$ is sharp in many situations, depending on the summability of the reaction $f$. Very recently, the range of attainable H\"older exponents has been widened in \cite{BDLMBS} for the homogeneous degenerate equation ($p>2$, $f=0$):
\begin{enumroman}
\item if $s\le (p-2)/p$, then $u\in C^\alpha_{\rm loc}(\Omega)$ for all $0<\alpha<ps/(p-2)$;
\item if $s>(p-2)/p$, then $u\in C^\alpha_{\rm loc}(\Omega)$ for all $0<\alpha<1$ (almost-Lipschitz continuity).
\end{enumroman}
In most results of this type, higher H\"older continuity is obtained through higher Sobolev regularity, see \cite{BL} for the degenerate case and \cite{DKLN} for the singular case. Being based on discrete differentiation of $u$ and integrability of $\nabla u$, such results are not immediately extended to operators with more general kernels (we will come back to this at the end of Section \ref{sec3}). The picture presented above is introductory and far from being complete, for instance see \cite{CS1,ROS1} for fully nonlinear nonlocal equations and \cite{L3} for the viscosity approach to the fractional $p$-Laplacian.

\section{Boundary regularity}\label{sec3}

\noindent
Following more closely the line of the present study, we want to examine the regularity of the solution up to the boundary $\partial\Omega$, which is assumed to be smooth. Already in the linear case, simple examples show that the regularity of the solution of fractional equations lowers dramatically as it approaches the boundary. To fix ideas, first consider the fractional Laplace equation on the half-line, coupled with a homogeneous nonlocal Dirichlet condition:
\[\begin{cases}
(-\Delta)^su = 0 & \text{in $(0,\infty)$} \\
u = 0 & \text{in $(-\infty,0]$.}
\end{cases}\]
The positive solution $u(x)=x_+^s$ is of class $C^\infty$ in $(0,\infty)$, but at best $C^s$ at $x=0$. Another well-known example is the torsion equation in a ball:
\[\begin{cases}
(-\Delta)^s u = 1 & \text{in $B_1(0)$} \\
u = 0 & \text{in $B_1^c(0)$.}
\end{cases}\]
The solution $u(x)=C_{N,s}(1-|x|^2)_+^s$ (see \cite{ROS} for the precise value of $C_{N,s}>0$) again is no more than $s$-H\"older continuous in $\bar B_1(0)$. So, the interior estimates of Theorem \ref{irl} are not stable at the boundary. In fact, it can be proved that $s$-H\"older in $\overline\Omega$ (and hence in $\R^N$) is the optimal regularity for the solutions of Dirichlet problems as soon as $\Omega$ satisfies the exterior ball condition, see \cite{ROS} for the fractional Laplacian and \cite{RO,ROS2} for general translation-invariant linear operators. In the proofs of boundary regularity results, a fractional Kelvin transform is often used.
\vskip2pt
\noindent
Turning to the nonlinear framework, no explicit solutions are known. Nevertheless, let us look at the estimate of Theorem \ref{irn}, taking for simplicity $\alpha=ps/(p-1)$. As soon as $x$ approaches $\partial\Omega$, we have $R\to 0$. The second term of the right-hand side then reduces to $\|f\|_{L^\infty(\Omega)}$, but the first and third term may blow up, depending on the behavior of $u$. Note that the nonlinear nature of the operator prevents use of transforms.
\vskip2pt
\noindent
The first global regularity result for the fractional $p$-Laplacian goes back to \cite{IMS}, and it ensures that the (unique) solution of the Dirichlet problem
\beq\label{dpp}
\begin{cases}
(-\Delta)_p^s u = f(x) & \text{in $\Omega$} \\
u = 0 & \text{in $\Omega^c$,}
\end{cases}
\eeq
with $\Omega$ bounded and $C^{1,1}$-smooth and $f\in L^\infty(\Omega)$, satisfies $u\in C^\alpha(\overline\Omega)$ for some indetermined $\alpha\in(0,1)$ depending on $N,p,s,\Omega,$ and $f$. Despite this result has been overcome by the following literature, one intermediate lemma \cite[Theorem 4.4]{IMS} has proved to be useful:

\begin{lemma}\label{bes}
Let $\Omega$ be bounded and $C^{1,1}$-smooth, $f\in L^\infty(\Omega)$, $u\in W^{s,p}_0(\Omega)$ be the solution of \eqref{dpp}, and set for all $x\in\R^N$
\[{\rm d}_\Omega(x) = {\rm dist}(x,\Omega^c).\]
Then, there exists $C=C(N,p,s,\Omega)>0$ s.t.\ for a.e.\ $x\in\Omega$
\[|u(x)| \le C\|f\|_{L^\infty(\Omega)}^\frac{1}{p-1}{\rm d}_\Omega^s(x).\]
\end{lemma}

\noindent
The proof of Lemma \ref{bes} divides in two steps: first, starting from the solution of a torsion problem on a ball we produce a barrier $w$ near a boundary point $x\in\partial\Omega$, s.t.\ $w\sim{\rm d}_\Omega^s$ near $x$; then, we apply a weak comparison principle to (a normalized) $u$ and such $w$. The estimate of Lemma \ref{bes} becomes meaningful near $\partial\Omega$, where it can be used along with Theorem \ref{irn} to prove the following global regularity result for the fractional $p$-Laplacian \cite[Theorem 2.7]{IM}:

\begin{theorem}\label{brn}
{\rm (Global regularity)} Let $\Omega$ be bounded and $C^{1,1}$-smooth, $f\in L^\infty(\Omega)$, $u\in W^{s,p}_0(\Omega)$ be the solution of \eqref{dpp}. Then, $u\in C^s(\overline\Omega)$ and there exists $C=C(N,p,s,\Omega)>0$ s.t.\
\[\|u\|_{C^s(\overline\Omega)} \le C\|f\|_{L^\infty(\Omega)}^\frac{1}{p-1}.\]
\end{theorem}

\noindent
The proof of Theorem \ref{brn} exploits Theorem \ref{irn} and Lemma \ref{bes} above: assuming for simplicity $\|f\|_{L^\infty(\Omega)}=1$, by Lemma \ref{bes} we have $\|u\|_{L^\infty(\Omega)}\le C{\rm diam}(\Omega)^s$. Also, given $x\in\partial\Omega$ and $R>0$ the same estimate can be used to see that
\[\underset{B_R(x)\cap\Omega}{\rm osc}\,u \le CR^s.\]
Now we invoke Theorem \ref{irn} with $\alpha=s$, so for all $B_{4R}(x)\subset\Omega$ we get $u\in C^s(B_{R/8}(x))$ and, using the estimates above, we find
\begin{align*}
[u]_{C^s(B_{R/8}(x))} &\le \frac{C}{R^s}\big[\|u\|_{L^\infty(B_R(x))}+R^\frac{ps}{p-1}+{\rm Tail}(u,x,R)\big] \\
&\le C+CR^\frac{s}{p-1}+{\rm diam}(\Omega)^\frac{s}{p-1}\Big[\int_{\Omega\cap B_R^c(x)}\frac{|u(y)|^{p-1}}{|x-y|^{N+ps}}\,dy\Big]^\frac{1}{p-1} \le C.
\end{align*}
So, the $C^s$-estimate is stable as $x$ approaches $\partial\Omega$. These conditions, via a technical lemma on H\"older continuous functions \cite[Lemma 2.6]{IM}, imply $u\in C^s(\overline\Omega)$ and the uniform estimate.
\vskip2pt
\noindent
Theorem \ref{brn} is optimal, even for the linear case $p=2$, due to the previous examples. Nevertheless, global $C^s$-regularity is not satisfactory, and in the absence of a general gradient bound, we are led to the study of a form or {\em fine} (or weighted) boundary regularity, which amounts at
\[\frac{u}{{\rm d}_\Omega^s} \in C^\alpha(\overline\Omega) \ \text{for some $\alpha\in(0,1)$,}\]
meaning of course that $u/{\rm d}_\Omega^s$ admits a $\alpha$-H\"older continuous extension to $\overline\Omega$. Such type of estimate comes from fully nonlinear regularity theory, and was first introduced in \cite{K} (with $s=1$) to prove $C^2$-regularity of solutions of second order elliptic equations. In the nonlocal framework, it was first proved in \cite{ROS} for $(-\Delta)^s$, in \cite{ROS2} for translation invariant linear operators, in \cite{KW} for linear operators with H\"older continuous kernels, and in \cite{ROS1} for fully nonlinear fractional equations. See also \cite{G} for nearly-optimal boundary smoothness for linear stable operators.
\vskip2pt
\noindent
For the fractional $p$-Laplacian, the path to fine boundary regularity is a long one. This result was achieved in \cite{IMS1} for the degenerate case $p>2$, and in \cite{IM} for the singular case $p<2$:

\begin{theorem}\label{fbr}
{\rm (Fine boundary regularity)} Let $\Omega$ be bounded and $C^{1,1}$-smooth, $f\in L^\infty(\Omega)$, $u\in W^{s,p}_0(\Omega)$ be the solution of \eqref{dpp}. Then, there exist $\alpha\in(0,s)$ and $C>0$ depending on $N,p,s,\Omega$, s.t.\ $u/{\rm d}_\Omega^s\in C^\alpha(\overline\Omega)$ and
\[\Big\|\frac{u}{{\rm d}_\Omega^s}\Big\|_{C^\alpha(\overline\Omega)} \le C\|f\|_{L^\infty(\Omega)}^\frac{1}{p-1}.\]
\end{theorem}

\noindent
Considering the lack of  boundary regularity inherent to nonlocal problems, Theorem \ref{fbr} represent the proper analogue to the celebrated $C^{1,\alpha}$-regularity result for the $p$-Laplacian from \cite{L2}. The proof is quite technical, but we will try to summarize it. We only consider the singular case $p<2$, which is more involved, essentially because the map $t\mapsto|t|^{p-2}t$ is non-differentiable at $0$. First we recall a typically nonlocal {\em superposition principle}:

\begin{lemma}\label{nsp}
Let $u\in\widetilde{W}^{s,p}(\Omega)$, $v\in L^1_{\rm loc}(\R^N)$ s.t.\ $\Omega\Subset V={\rm supp}(u-v)$ and
\[\int_V \frac{|v(y)|^{p-1}}{(1+|y|)^{N+ps}}\,dy < \infty.\]
Set for all $x\in\R^N$
\[w(x) = \begin{cases}
u(x) & \text{if $x\in V^c$} \\
v(x) & \text{if $x\in V$.}
\end{cases}\]
Then, $w\in\widetilde{W}^{s,p}(\Omega)$ and weakly for $x\in\Omega$
\[(-\Delta)_p^s w(x) = (-\Delta)_p^s u(x)+2\int_V \frac{(u(x)-v(y))^{p-1}-(u(x)-u(y))^{p-1}}{|x-y|^{N+ps}}\,dx\footnote{For brevity we set $a^{p-1}=|a|^{p-2}a$.}.\]
\end{lemma}

\noindent
The main idea is to perturb (super, sub) solutions at a distance from the domain, and obtain through Lemma \ref{nsp} a control on the fractional $p$-Laplacians of the perturbed functions inside the domain. More precisely, assume for simplicity that $0\in\partial\Omega$ and set $D_R=B_R(0)\cap\Omega$, while $\tilde B_R$ denotes a ball of radius $R/4$ inside $\Omega$, centered along the normal direction to $\partial\Omega$, s.t.\ ${\rm d}_\Omega$ is controlled by multiples of $R$ in $\tilde B_R$. We aim at an oscillation estimate for the quotient $v=u/{\rm d}_\Omega^s$ in $D_R$. We define the {\em nonlocal excess} with respect to a constant $m\in\R$ as
\[L(u,m,R) = \Big[\dashint_{\tilde B_R}|v(x)-m|^{p-1}\,dx\Big]^\frac{1}{p-1}.\]
The quantity above controls $v$ in $D_{R/2}$ by means of two weak Harnack inequalities. Namely, let $u\in\widetilde{W}^{s,p}(D_R)$ satisfy for some $m,K,H>0$
\[\begin{cases}
(-\Delta)_p^s u \ge -\min\{K,H\} & \text{in $D_R$} \\
u \ge m{\rm d}_\Omega^s & \text{in $\R^N$}
\end{cases}\]
(note the {\em global} lower bound requested for $u$). Then, there exist $\sigma\in(0,1)$, $C>0$ only depending on $N,p,s,\Omega$ s.t.\
\[\inf_{D_{R/2}}(v-m) \ge \sigma L(u,m,R)-C(KR^s)^\frac{1}{p-1}-C(m+Hm^{2-p})R^s.\]
This lower bound is obtained by superposition between $u$ and convenient barrier functions, respectively, a torsion solution when the excess is small, and a diffeomorphic equivalent to the distance from $\partial\Omega$ (this is where the $C^{1,1}$-smoothness assumptions on $\partial\Omega$ is mainly required).
\vskip2pt
\noindent
An apparently equivalent upper bound holds for subsolutions. If $u$ satisfies
\[\begin{cases}
(-\Delta)_p^s u \le \min\{K,H\} & \text{in $D_R$} \\
u \le M{\rm d}_\Omega^s & \text{in $\R^N$,}
\end{cases}\]
then for some universal constants we have
\[\inf_{D_{R/4}}(M-v) \ge \sigma L(u,M,R)-C(KR^s)^\frac{1}{p-1}-C(M+HM^{2-p})R^s.\]
The main difference is that the bound $u\le M{\rm d}_\Omega^s$ leaves the sign of $u$ undetermined. We use the fact that, if the excess is big enough, then $u\le 0$ in $D_{R/2}$. We are now able to prove the desired oscillation estimate on $v$, in the following form: define a sequence of radii $(R_n)$ tending to $0$, then there exist two more sequences, $(m_n)$ nondecreasing and $(M_n)$ nonincreasing, s.t.\ for all $n\in\N$
\[m_n \le \inf_{D_{R_n}}v \le \sup_{D_{R_n}}v \le M_n, \ M_n-m_n = \mu R_n^\alpha,\]
with $\alpha\in(0,1)$, $\mu>1$ only depending on the data. The above estimates are proved by strong induction. For $n=0$, it is simply Lemma \ref{bes} above. Assuming that the relation above holds for all integers between $0$ and $n$, we apply the lower bound to $(u\vee m_n{\rm d}_\Omega^s)$ and the upper bound to $(u\wedge M_n{\rm d}_\Omega^s)$ (which satisfy global bounds), and, via delicate estimates on excess and tail terms, we find $m_{n+1}$, $M_{n+1}$. After that, it is standard to find $C>0$ s.t.\ for all $r>0$ small enough
\[\underset{D_r}{\rm osc}\,v \le C\|f\|_{L^\infty(\Omega)}^\frac{1}{p-1}r^\alpha,\]
which in turn implies $v\in C^\alpha(\overline\Omega)$ as in Theorem \ref{brn}, thus concluding the proof of Theorem \ref{fbr}.
\vskip2pt
\noindent
In \cite{MS}, a stimulating list of open problems in fractional regularity theory was proposed. Since then, much has been achieved, but some limit cases still lack a solution, so we would like to provide the reader with an updated report. Problem $(1)$ from \cite{MS} deals with higher interior H\"older regularity for the solution of
\[(-\Delta)_p^s u = f(x) \ \text{in $\Omega$,}\]
with $f\in L^\infty(\Omega)$, proposing the conjecture that $u\in C^\alpha_{\rm loc}(\Omega)$ for some $\alpha>s$. We know from Theorem \ref{irn} that this is true. Further, the supremum of attainable H\"older exponents $ps/(p-1)$ is optimal for $s<(p-1)/p$, see \cite{BLS,GL}. On the contrary, for $s>(p-1)/p$ no more than Lipschitz regularity is known so far.
\vskip2pt
\noindent
Problem $(2)$ is about higher differentiability (or Sobolev regularity) of the solutions. The result of \cite{BL} has been extended and improved in several ways. The most general results we are aware of are in \cite{DKLN}, covering both the singular and the degenerate cases and general reactions, proving the following implications for any $p>1$, $f\in L^{p'}(\Omega)$:
\begin{enumroman}
\item if $p\ge 2$ or $p<2$ and $s<(p-1)/p$, then $u\in W^{\sigma,p}(\Omega)$ for all $\sigma\in(0,ps/(p-1))$;
\item if $p<2$ and $s>(p-1)/p$, then $u$ is differentiable with $\nabla u\in W^{\sigma,p}(\Omega)$ for $\sigma\in(0,ps-p+1)$.
\end{enumroman}
Even better regularity is obtained if $f\in W^{\tau,p'}(\Omega)$ for convenient $\tau$, see \cite[Theorems 1.5, 1.7]{DKLN}. Nevertheless, optimal conditions on $f$ are still missing.
\vskip2pt
\noindent
Problem $(3)$ tackles the boundary regularity of the solution of the Dirichlet problem \eqref{dpp} under the following general assumption on the smoothness of $\partial\Omega$: there exists $R>0$ s.t.\
\[\inf_{x\in\partial\Omega}\,\inf_{r\in(0,R)}\,\frac{|B_r(x)\cap\Omega^c|}{r^N} > 0.\]
The condition above was introduced in \cite{KKP}, where H\"older continuity up to the boundary is proved for the solution of an obstacle problem for the homogeneous fractional $p$-Laplace equation. To the best of our knowledge, it is still an open problem whether the obstacle can be replaced with a non-zero reaction $f\in L^\infty(\Omega)$.
\vskip2pt
\noindent
Finally, problem $(4)$ deals with fine boundary regularity and it has found a complete answer in Theorem \ref{fbr} above.
\vskip2pt
\noindent
Inspired by \cite{MS}, we would like to leave the reader with a new short list of open problems related to the main subject of this note, i.e., boundary regularity:

\begin{problem}\label{p1}
Find the optimal $\alpha$ in Theorem \ref{fbr}. For $p=2$, we know from \cite{ROS} that $u/{\rm d}_\Omega^s\in C^\alpha(\overline\Omega)$ for all $\alpha\in(0,s)$. We may reasonably conjecture that $s$ is the supremum of H\"older exponents in the nonlinear case as well, but a careful reading of the proofs in \cite{IM,IMS1} will show that the exponent $\alpha$ found there is far from being explicit.
\end{problem}

\begin{problem}\label{p2}
Prove boundary regularity for the fractional $p$-Laplacian with a reaction $f\in L^q(\Omega)$, $q>\max\{1,N/ps\}$. Note that, as soon as $f$ is unbounded, Lemma \ref{bes} fails. Clearly, the problem of {\em fine} boundary regularity for unbounded reactions would come next.
\end{problem}

\begin{problem}\label{p3}
Find conditions on the kernel $K$ s.t.\ fine boundary regularity holds for the general operator $\mathcal{L}_{p,K}^s$. In general, hypotheses $(K_1)$, $(K_2)$ are not enough even in the linear case, but in \cite{RO} it is proved that $u/{\rm d}_\Omega^s\in C^\alpha(\overline\Omega)$ for all $\alpha\in(0,s)$, provided the kernel is of the form
\[K(x,y) = \frac{1}{|x-y|^{N+2s}}\mu\Big(\frac{x-y}{|x-y|}\Big),\]
where $\mu$ is a measure on the sphere $\mathbb{S}^{N-1}$ (better regularity is achieved if $\mu$ is smoother). We do not know if this type of kernels have been treated also in the nonlinear framework.
\end{problem}

\section{Applications}\label{sec4}

\noindent
Just like the result of \cite{L2}, Theorem \ref{fbr} is especially suitable for applications in nonlinear analysis, meaning the variety of variational and topological methods to prove existence, multiplicity, and qualitative properties of the solutions of general problems of the form:
\beq\label{dpn}
\begin{cases}
(-\Delta)_p^s u = f(x,u) & \text{in $\Omega$} \\
u = 0 & \text{in $\Omega^c$,}
\end{cases}
\eeq
where $f:\Omega\times\R\to\R$ is a Carath\'eodory mapping subject to certain growth and/or monotonicity assumptions (see \cite{MMP} for a comprehensive introduction to such methods for nonlinear elliptic equations). For simplicity, we assume that $f(x,\cdot)$ has {\em subcritical growth}, i.e., for a.e.\ $x\in\Omega$ and all $t\in\R$
\beq\label{gc}
|f(x,t)| \le C_0(1+|t|^{r-1}) \ (C_0>0,\,r\in(1,p^*_s)).
\eeq
We derive from \cite[Proposition 2.3]{IM1} a uniform {\em a priori} bound for the solutions of problem \eqref{dpn} (which are defined as in Section \ref{sec1}):

\begin{proposition}\label{apb}
{\rm (A priori bound)} Let $f$ satisfy \eqref{gc}, $u\in W^{s,p}_0(\Omega)$ be a solution of \eqref{dpn}. Then, $u\in L^\infty(\Omega)$ and there exists $C=C(N,p,s,\Omega,r,C_0,\|u\|_{W^{s,p}_0(\Omega)})>0$ s.t.\ $\|u\|_{L^\infty(\Omega)}\le C$.
\end{proposition}

\noindent
By \eqref{gc} and Proposition \ref{apb} above, it is easily seen that $f(\cdot,u)\in L^\infty(\Omega)$. Therefore, by Theorem \ref{fbr} we have $u/{\rm d}_\Omega^s\in C^\alpha(\overline\Omega)$, with $\alpha\in(0,1)$ depending on the data and a uniform bound on the $C^\alpha(\overline\Omega)$-norm of the quotient. We may conveniently rephrase such result by saying that $u\in C^\alpha_s(\overline\Omega)$, where we have set
\[C^\alpha_s(\overline\Omega) = \Big\{u\in C^0(\overline\Omega):\,\frac{u}{{\rm d}_\Omega^s}\in C^\alpha(\overline\Omega)\Big\},\]
endowed with the norm $\|u\|_{C^\alpha_s(\overline\Omega)}=\|u/{\rm d}_\Omega^s\|_{C^\alpha(\overline\Omega)}$. This weighted H\"older space plays, in fractional elliptic boundary value problems, the same role as $C^{1,\alpha}(\overline\Omega)$ in local $p$-Laplacian problems. Note that for all $\alpha\in(0,1)$ $C^\alpha_s(\overline\Omega)$ is compactly embedded into $C^0_s(\overline\Omega)$, and the latter has a positive order cone with nonempty interior
\[{\rm int}(C^0_s(\overline\Omega)_+) = \Big\{u\in C^0_s(\overline\Omega):\,\inf_{\Omega}\frac{u}{{\rm d}_\Omega^s}>0\Big\}.\]
This functional framework bears several applications in the study of problems of the class \eqref{dpn}. First, we recall from \cite[Theorems 2.6, 2.7]{IMP} the following strong mimimum/comparison principles (where inequalities are meant in the weak sense):

\begin{theorem}\label{mcp}
Let $f\in C^0(\R)\cap BV_{\rm loc}(\R)$. Then we have:
\begin{enumroman}
\item\label{mcp1} {\rm (Strong minimum principle)} if $u\in\widetilde{W}^{s,p}(\Omega)\cap C^0(\overline\Omega)$, $u\neq 0$ satisfies
\[\begin{cases}
(-\Delta)_p^s u+f(u) \ge f(0) & \text{in $\Omega$} \\
u \ge 0 & \text{in $\Omega^c$,}
\end{cases}\]
then
\[\inf_\Omega\,\frac{u}{{\rm d}_\Omega^s} > 0;\]
\item\label{mcp2} {\rm (Strong comparison principle}) if $u,v\in\widetilde{W}^{s,p}(\Omega)\cap C^0(\overline\Omega)$, $u\neq v$, $C>0$ satisfy
\[\begin{cases}
(-\Delta)_p^s v+f(v) \le (-\Delta)_p^s u+f(u) \le C & \text{in $\Omega$} \\
0 = v \le u & \text{in $\Omega^c$,}
\end{cases}\]
then
\[\inf_\Omega\,\frac{u-v}{{\rm d}_\Omega^s} > 0.\]
\end{enumroman}
\end{theorem}

\noindent
The proof of Theorem \ref{mcp} relies on the superposition principle of Lemma \ref{nsp}. As a consequence of nonlocal diffusion of positivity, the minimum/comparison principles above yield a more detailed information about the behavior of (super) solutions, with respect to previous results for fractional operators, see for instance \cite{DPQ,J}. Incidentally we note that the strong minimum principle \ref{mcp1} even holds for equations like
\[(-\Delta)_p^s u = u^{q-1} \ (1<q<p),\]
preventing the appearance of {\em dead cores} without requiring additional integrability conditions, unlike in the case of the $p$-Laplacian (see \cite{V}). Similarly, the strong comparison principle \ref{mcp2} gives a global information unavailable in the local case (see \cite{GV}).
\vskip2pt
\noindent
When referred to {\em solutions}, Theorem \ref{mcp} \ref{mcp1} rephrases in the words of the functional framework introduced above, ensuring under suitable conditions that any non-negative, non-zero solution $u\in W^{s,p}_0(\Omega)$ of \eqref{dpn} satisfies $u \in {\rm int}(C^0_s(\overline\Omega)_+)$. This can be seen as a fractional version of Hopf's boundary point lemma, i.e., for all $x\in\partial\Omega$
\[\frac{\partial u}{\partial\nu^s}(x) \sim \lim_{y\to x}\frac{u(y)}{{\rm d}_\Omega^s(y)} > 0.\]
Problem \eqref{dpn} can be studied following essentially two approaches. In the {\em variational} approach, we seek solutions as critical points of an energy functional (for semilinear fractional equations, we refer to \cite{MBRS}). Set for all $(x,t)\in\Omega\times\R$
\[F(x,t) = \int_0^t f(x,\tau)\,d\tau,\]
and for all $u\in W^{s,p}_0(\Omega)$
\[\Phi(u) = \frac{\|u\|_{W^{s,p}_0(\Omega)}^p}{p}-\int_\Omega F(x,u)\,dx.\]
Then, $\Phi\in C^1(W^{s,p}_0(\Omega))$ and the solutions of \eqref{dpn} coincide with the critical points of $\Phi$, in particular its local minimizers.
\vskip2pt
\noindent
In many situations, one is lead to deal with truncations of $f(x,\cdot)$ and the corresponding truncated functionals, which coincide with $\Phi$ only on certain order-related subsets of $W^{s,p}_0(\Omega)$ (a typical example is the positive order cone). Due to the nature of the $W^{s,p}_0(\Omega)$-topology, such sets usually have empty interior. Thus, it is useful to swap to some H\"older type space in order to minimize truncated functionals on nonempty open sets (this method was introduced in \cite{BN} to deal with problems at critical growth), preserving the local minimizers.
\vskip2pt
\noindent
Such change of topology is made possible by the coincidence of local minimizers in $W^{s,p}_0(\Omega)$ and in $C^0_s(\overline\Omega)$, respectively, which was proved in general in \cite{IM1}:

\begin{theorem}\label{svh}
{\rm (Sobolev vs.\ H\"older minima)} Let $f$ satisfy \eqref{gc}, $u\in W^{s,p}_0(\Omega)$. Then, the following are equivalent:
\begin{enumroman}
\item\label{svh1} there exists $\rho>0$ s.t.\ $\Phi(u+v)\ge\Phi(u)$ for all $v\in W^{s,p}_0(\Omega)$, $\|v\|_{W^{s,p}_0(\Omega)}\le\rho$;
\item\label{svh2} there exists $\sigma>0$ s.t.\ $\Phi(u+v)\ge\Phi(u)$ for all $v\in W^{s,p}_0(\Omega)\cap C^0_s(\overline\Omega)$, $\|v\|_{C^0_s(\overline\Omega)}\le\sigma$.
\end{enumroman}
\end{theorem}

\noindent
The proof of Theorem \ref{svh} relies on the $C^\alpha_s(\overline\Omega)$-bound of Theorem \ref{fbr} and some uniform {\em a priori} bounds on perturbed problems, and it even holds for reactions at critical growth.
\vskip2pt
\noindent
Exploiting Theorems \ref{mcp} and \ref{svh}, we can prove existence and multiplicity of solutions under a variety of assumptions, see for instance \cite[Theorem 5.3]{ILPS} (where fine boundary regularity was {\em assumed} in a slightly more general form, while at present it is ensured by Theorem \ref{fbr}). A more recent example is the following, dealing with $(p-1)$-superlinear reactions at infinity, see \cite[Theorem 1.1]{ISV}:

\begin{theorem}\label{sup}
Let $f$ satisfy \eqref{gc}, and in addition uniformly for a.e.\ $x\in\Omega$
\begin{enumroman}
\item\label{sup1} $\displaystyle\lim_{|t|\to\infty}\frac{F(x,t)}{|t|^p} = \infty$;
\item\label{sup2} $\displaystyle\liminf_{|t|\to\infty}\frac{f(x,t)t-pF(x,t)}{|t|^q} > 0$ for some $\displaystyle q\in\Big(\frac{N(r-p)}{ps},p^*_s\Big)$;
\item\label{sup3} $\displaystyle\lim_{t\to 0}\frac{f(x,t)}{|t|^{p-2}t} = 0$.
\end{enumroman}
Then, problem \eqref{dpn} has at least three nontrivial solutions $u_+>0$, $u_-<0$, and $\tilde u\neq 0$.
\end{theorem}

\noindent
To prove Theorem \ref{sup}, we truncate the reaction at the origin, thus introducing
\[f_\pm(x,t) = f(x,\pm t^\pm),\]
which still satisfies \eqref{gc} and one-sided versions of the hypotheses \ref{sup1}-\ref{sup3}. The corresponding energy functionals are denoted $\Phi_\pm$. Each has a strict local minimizer at $0$ and a mountain pass critical point $u_\pm\in\pm{\rm int}(C^0_s(\overline\Omega)_+)$ (Theorem \ref{mcp}), and by Theorem \ref{svh} the same situation reflects on $\Phi$. Finally, a Morse-theoretic argument reveals the existence of a fourth critical point $\tilde u$.
\vskip2pt
\noindent
The {\em topological} approach to problem \eqref{dpn} is essentially based on the degree theory for $(S)_+$-operators, see \cite{MMP}. We define a nonlinear operator $A:W^{s,p}_0(\Omega)\to W^{-s,p'}(\Omega)$ by setting for all $u,\varphi\in W^{s,p}_0(\Omega)$
\[\langle A(u),\varphi\rangle = \iint_{\R^N\times\R^N}\frac{|u(x)-u(y)|^{p-2}(u(x)-u(y))(\varphi(x)-\varphi(y))}{|x-y|^{N+ps}}\,dx\,dy-\int_\Omega f(x,u)\varphi\,dx,\]
and we seek zeros of $A$ in $W^{s,p}_0(\Omega)$ (the two approaches are in principle equivalent, as $A=\Phi'$). Degree theory is used, for instance, in \cite{I1}, where several existence results are presented under conditions relating the asymptotic behavior of $f(x,\cdot)$ to the principal eigenvalue $\lambda_1$ of the fractional $p$-Laplacian in $W^{s,p}_0(\Omega)$ (see for instance \cite{FP,I,LL} on the variational spectrum of $(-\Delta)_p^s$). The following is a slightly simplified version of \cite[Theorem 5.2]{I1} for asymptotically $(p-1)$-linear reactions:

\begin{theorem}\label{deg}
Let $f$ satisfy \eqref{gc}, and in addition uniformly for a.e.\ $x\in\Omega$
\begin{enumroman}
\item\label{deg1} $\displaystyle\lim_{t\to\infty}\frac{f(x,t)}{t^{p-1}} \in (\lambda_1,\infty)$;
\item\label{deg2} $\displaystyle\lim_{t\to 0^+}\frac{f(x,t)}{t^{p-1}} \in (\lambda_1,\infty)$;
\item\label{deg3} $\displaystyle\limsup_{t\to c^-}\frac{f(x,t)}{(c-t)^{p-1}} \le K$ for some $c,K>0$.
\end{enumroman}
Then, problem \eqref{dpn} has at least two positive solutions $u_1,u_2>0$.
\end{theorem}

\noindent
In the proof of Theorem \ref{deg}, the use of the constant $c$ as a supersolution plays a crucial role, as well as the strong comparison principle of Theorem \ref{mcp} \ref{mcp2}. The first solution $u_1$ is found by minimizing a conveniently truncated energy functional, then the second solution is achieved via a topological argument by contradiction, by computing the degree of the operator $A$ in a large ball $B_R(0)$ and in two small balls $B_r(0)$, $B_r(u_1)$, respectively.
\vskip2pt
\noindent
Another interesting application of our regularity results is connected with the study of {\em extremal solutions} of problem \eqref{dpn} in a sub-supersolution interval. First, let $\underline u,\overline u\in\widetilde{W}^{s,p}(\Omega)$ be a sub- and a supersolution of \eqref{dpn}, respectively, s.t.\ $\underline u\le\overline u$ in $\Omega$. The solution set
\[\mathcal{S}(\underline u,\overline u) = \big\{u\in W^{s,p}_0(\Omega):\,\text{$u$ solution of \eqref{dpn}, $\underline u\le u\le\overline u$ in $\Omega$}\big\}\]
is nonempty, due to classical results in operator theory (see \cite{MMP}). In \cite{FI}, the topological properties of such set are investigated:

\begin{theorem}\label{ssi}
{\rm (Solutions in a sub-supersolution interval)} Let $f$ satisfy \eqref{gc}, $\underline u,\overline u\in\widetilde{W}^{s,p}(\Omega)$ be a sub- and a supersolution of \eqref{dpn} s.t.\ $\underline u\le\overline u$ in $\Omega$. Then, $\mathcal{S}(\underline u,\overline u)$ has the following properties:
\begin{enumroman}
\item\label{ssi1} $\mathcal{S}(\underline u,\overline u)$ is both upward and downward directed;
\item\label{ssi2} $\mathcal{S}(\underline u,\overline u)$ is compact in both $W^{s,p}(\Omega)$ and $C^0_s(\overline\Omega)$;
\item\label{ssi3} there exist $u_1,u_2\in\mathcal{S}(\underline u,\overline u)$ s.t.\ $u_1\le u\le u_2$ in $\Omega$ for all $u\in\mathcal{S}(\underline u,\overline u)$.
\end{enumroman}
\end{theorem}

\noindent
The proof of Theorem \ref{ssi} is long but straightforward. First, we endow $W^{s,p}_0(\Omega)$ with a lattice structure, see \cite{GM}, and we prove that whenever $u,v$ are subsolutions (resp., supersolutions) of \eqref{dpn}, then $u\vee v$ is as well a subsolution (resp., $u\wedge v$ is a supersolution). This implies \ref{ssi1}. To prove \ref{ssi2}, we use the uniform bound in $C^\alpha_s(\overline\Omega)$ from Theorem \ref{fbr} and the compact embedding of $C^\alpha_s(\overline\Omega)$ into $C^0_s(\overline\Omega)$ to prove compactness in the latter space, then the $(S)_+$-property of the fractional $p$-Laplacian to achieve compactness in $W^{s,p}_0(\Omega)$ as well. Finally, there are several possible arguments, based on either Cantor's diagonal ordering or Zorn's lemma, to prove the existence of extremal elements \ref{ssi3}.
\vskip2pt
\noindent
Theorem \ref{ssi} can be employed in many ways to prove the existence of extremal solutions of fractional $p$-Laplacian equations in a given sub-supersolution interval, see \cite{IM1,IMP} for some applications to logistic equations. In particular, we want to recall here a result ensuring existence of the smallest positive solution (this is a one-sided version of \cite[Theorem 4.1]{FI}):

\begin{theorem}\label{ext}
Let $f$ satisfy \eqref{gc}, and in addition uniformly for a.e.\ $x\in\Omega$
\begin{enumroman}
\item\label{ext1} $\displaystyle\limsup_{t\to \infty}\frac{F(x,t)}{t^p} < \frac{\lambda_1}{p}$;
\item\label{ext2} $\displaystyle\lim_{t\to 0^+}\frac{f(x,t)}{t^{p-1}} \in (\lambda_1,\infty)$.
\end{enumroman}
Then, problem \eqref{dpn} has a smallest positive solution $u_+>0$.
\end{theorem}

\noindent
To prove Theorem \ref{ext}, once again we begin by minimizing a truncation of the energy functional $\Phi$, which via Theorems \ref{mcp} and \ref{svh} leads to the existence of a positive solution $u_1\in{\rm int}(C^0_s(\overline\Omega)_+)$. Combining such $u_1$ with the positive principal eigenfunction of $(-\Delta)_p^s$ we find a sequence of positive subsolutions and hence, by Theorem \ref{ssi}, a decreasing sequence of minimal solutions $(u_n)$. A compactness argument then allows to make a subsequence of $(u_n)$ converge to $u_+$, which turns out to be the smallest positive solution of \eqref{dpn}.
\vskip2pt
\noindent
Under symmetric hypotheses, it is possible to prove in a similar way the existence of the biggest negative solution $u_-<0$, and then via the mountain pass theorem, of a nodal solution $\tilde u\neq 0$ s.t.\ $u_-\le\tilde u\le u_+$ in $\Omega$ (see \cite[Theorem 5.1]{FI}).

\begin{remark}\label{sin}
We note that, in the cited references, Theorems \ref{deg}, \ref{ssi}, and \ref{ext} were only stated for $p\ge 2$. This is due to regularity issues only, since fine boundary regularity was only available for the degenerate fractional $p$-Laplacian at that time. Now, the use of Theorem \ref{fbr} allows to prove such results in the singular case $p<2$ as well.
\end{remark}

\vskip4pt
\noindent
{\bf Acknowledgement.} The present note is based on a talk given by the author at the University of Bologna in February, 2024, in the series {\em Seminari di Analisi Matematica Bruno Pini}. The author is a member of GNAMPA (Gruppo Nazionale per l'Analisi Matematica, la Probabilit\`a e le loro Applicazioni) of INdAM (Istituto Nazionale di Alta Matematica 'Francesco Severi'), and is partially supported by the research project {\em Problemi non locali di tipo stazionario ed evolutivo} (GNAMPA, CUP E53C23001670001). The author wishes to thank the anonymous Referee for valuable comments and for drawing his attention to the interesting reference \cite{DTGCV}.

\end{document}